\documentclass [12pt, twoside]{amsart}
\usepackage{amsmath,amsfonts,amsthm,amssymb}

\newtheorem{theorem}{Theorem}[section]
\newtheorem{lemma}[theorem]{Lemma}
\newtheorem{proposition}[theorem]{Proposition}
\newtheorem{corollary}[theorem]{Corollary}
\newtheorem{question}[theorem]{Question}

\theoremstyle{definition}

\newtheorem{example}[theorem]{Example}

\theoremstyle{remark}
\newtheorem{remark}[theorem]{Remark}

\numberwithin{equation}{section}

\begin{document}

\title [Dynamics of tuples of matrices]{Dynamics of tuples of matrices}

\author[G. Costakis]{G. Costakis}
\address{Department of Mathematics, University of Crete, Knossos Avenue, GR-714 09 Heraklion, Crete, Greece}
\email{costakis@math.uoc.gr}
\thanks{}

\author[D. Hadjiloucas]{D. Hadjiloucas}
\address{The School of Sciences, European University Cyprus, 6 Diogenes Street, Engomi, P.O.Box 22006, 1516 Nicosia, Cyprus}
\email{d.hadjiloucas@euc.ac.cy}
\thanks{}

\author[A. Manoussos]{A. Manoussos}
\address{Fakult\"{a}t f\"{u}r Mathematik, SFB 701, Universit\"{a}t Bielefeld, Postfach 100131, D-33501 Bielefeld, Germany}
\email{amanouss@math.uni-bielefeld.de}
\thanks{During this research the third author was fully supported by SFB 701 ``Spektrale Strukturen und
Topologische Methoden in der Mathematik" at the University of Bielefeld, Germany. He would also like to express his gratitude to Professor H. Abels for his support.}



\begin{abstract}
In this article we answer a question raised by N. Feldman in \cite{Feldman} concerning the dynamics of tuples of operators on $\mathbb{R}^n$. In particular, we prove
that for every positive integer $n\geq 2$ there exist $n$ tuples $(A_1, A_2, \ldots ,A_n)$ of $n\times n$ matrices over $\mathbb{R}$ such that $(A_1, A_2, \ldots
,A_n)$ is hypercyclic. We also establish related results for tuples of $2\times 2$ matrices over $\mathbb{R}$ or $\mathbb{C}$ being in Jordan form.
\end{abstract}
\maketitle
\section{Introduction} \label{S1}
Following the recent work of Feldman in \cite{Feldman}, an $n$-tuple of operators is a finite sequence of length $n$ of commuting
continuous linear operators $T_1, T_2,\cdots, T_n$ acting on a locally convex space $X$. The tuple $(T_1, T_2, \ldots ,T_n)$ is
hypercyclic if there exists a vector $x\in X$ such that the set
$$\{ T_1^{k_1}T_2^{k_2}\cdots T_n^{k_n}x : k_1, k_2,\ldots k_n \geq 0 \}$$ is dense in $X$. Such a vector $x$ is
called hypercyclic for $(T_1,T_2, \ldots ,T_n)$ and the set of hypercyclic vectors for $(T_1,T_2, \ldots ,T_n)$
will be denoted by $HC((T_1,T_2, \ldots ,T_n))$. The above definition generalizes the notion of
hypercyclicity to tuples of operators. For an account of results, comments and an
extensive bibliography on hypercyclicity we refer to \cite{BM}, \cite{GE}, \cite{GE2} and \cite{GE3}. For
results concerning the dynamics of tuples of operators see \cite{Feldman1}, \cite{Feldman2}, \cite{Feldman} and
\cite{Ker}.

In \cite{Feldman} Feldman showed, among other things, that in $\mathbb{C}^n$ there exist diagonalizable
$n+1$-tuples of matrices having dense orbits. In addition he proved that there is no $n$-tuple of diagonalizable
matrices on $\mathbb{R}^n$ or $\mathbb{C}^n$ that has a somewhere dense orbit. Therefore the following question
arose naturally.\\
{\bf Question (Feldman \cite{Feldman})}. \textit{Are there non-diagonalizable $n$-tuples on $\mathbb{R}^k$ that
have somewhere dense orbits}?

We give a positive answer to this question in a very strong form, as the next theorem shows.

\begin{theorem} \label{maintheorem}
For every positive integer $n\geq 2$ there exist $n$-tuples $(A_1, A_2, \ldots ,A_n)$ of $n\times n$ non-diagonalizable matrices over $\mathbb{R}$ such that $(A_1,
A_2, \ldots ,A_n)$ is hypercyclic.
\end{theorem}

Restricting ourselves to tuples of $2 \times 2$ matrices in Jordan form either on $\mathbb{R}^2$
or $\mathbb{C}^2$, we prove the following.

\begin{theorem}\label{R2hypjordan}
There exist  $2\times 2$ matrices $A_j, j=1,2,3,4$ in Jordan form over $\mathbb{R}$ such that $(A_1, A_2, A_3,
A_4)$ is hypercyclic. In particular
\[
HC((A_1,A_2,A_3,A_4))= \left\{ \left(
\begin{array}{c}
x_1 \\
x_2
\end{array} \right) \in \mathbb{R}^2 : x_2\neq 0 \right\} .
\]
\end{theorem}

\begin{theorem}\label{C2hypjordan}
There exist  $2\times 2$ matrices $A_j, j=1,2,\ldots ,8$ in Jordan form over $\mathbb{C}$ such that $(A_1, A_2,
\ldots , A_8)$ is hypercyclic.
\end{theorem}

\section{Products of $2\times 2$ matrices}

\begin{lemma} \label{jordanproduct}
Let $m$ be a positive integer and for each $j=1,2,\ldots ,m$ let $A_j$ be a $2\times 2$ matrix in Jordan form over a field $\mathbb{F}=\mathbb{C} \,\, \textrm{or}
\,\, \mathbb{R}$ , i.e. $A_j=\left( \begin{array}{cc}
a_j & 1 \\
0 & a_j
\end{array} \right)$
for $a_1,a_2,\ldots,a_m\in\mathbb{F}$.
Then $(A_1, A_2,\ldots ,A_m)$ over $\mathbb{C}$ (respectively $\mathbb{R}$) is hypercyclic if and only if the
sequence

\[ \left\{\left( \begin{array}{c}
\frac{k_1}{a_1}+\frac{k_2}{a_2}+\ldots  +\frac{k_m}{a_m} \\
{a_1}^{k_1}{a_2}^{k_2}\ldots {a_m}^{k_m}
\end{array} \right) : k_1, k_2, \ldots ,k_m \in \mathbb{N}\right\} \]
is dense in $\mathbb{C}^2$ (respectively $\mathbb{R}^2$).
\end{lemma}
\begin{proof}
We prove the above in the case $\mathbb{F}=\mathbb{C}$, since the other case is similar. Observe that \[
{A_j}^l=\left(
\begin{array}{cc}
{a_j}^l & l {a_j}^{l-1}\\
0 & {a_j}^l
\end{array} \right)\] for $l\in \mathbb{N}$. As a result we have
\[ {A_1}^{k_1}{A_2}^{k_2}\ldots {A_m}^{k_m}=\left( \begin{array}{cc}
\prod\limits_{j=1}^m{a_j}^{k_j} &  \prod\limits_{j=1}^m{a_j}^{k_j}\sum\limits_{s=1}^{m}\frac{k_s}{a_s}\\
0 & \prod\limits_{j=1}^m{a_j}^{k_j} \end{array}  \right) .\]
Assume that $(A_1, A_2,\ldots ,A_m)$ is hypercyclic and let $ \left( \begin{array}{c}
z_1 \\
z_2
\end{array} \right) \in \mathbb{C}^2$ be a hypercyclic vector for $(A_1, A_2,\ldots ,A_m)$.
Then the sequence
\[ \left\{ {A_1}^{k_1}{A_2}^{k_2}\ldots {A_m}^{k_m} \left( \begin{array}{c}
z_1 \\
z_2
\end{array} \right): k_1, k_2, \ldots ,k_m \in \mathbb{N} \right\} \]
\[ = \left\{ \left( \begin{array}{c}
z_1\prod\limits_{j=1}^m{a_j}^{k_j} + z_2\prod\limits_{j=1}^m{a_j}^{k_j}\sum\limits_{s=1}^{m}\frac{k_s}{a_s}\\
z_2\prod\limits_{j=1}^m{a_j}^{k_j} \end{array}  \right) : k_1, k_2, \ldots ,k_m \in \mathbb{N} \right\} \]
is dense in $\mathbb{C}^2$. This implies that $z_2\neq 0$. Dividing the element on the first row by that on the
second, it can easily be shown that the sequence \[ \left\{\left( \begin{array}{c}
\frac{k_1}{a_1}+\frac{k_2}{a_2}+\ldots  +\frac{k_m}{a_m} \\
{a_1}^{k_1}{a_2}^{k_2}\ldots {a_m}^{k_m}
\end{array} \right) : k_1, k_2, \ldots ,k_m \in \mathbb{N}\right\} \]
is dense in $\mathbb{C}^2$. The converse can easily be shown.
\end{proof}

\begin{remark}\label{hypvectors}
Let $m$ be a positive integer and for each $j=1,2,\ldots ,m$ let $A_j$ be a $2\times 2$ matrix in Jordan form
over a field $\mathbb{F}=\mathbb{C} \,\, \textrm{or} \,\, \mathbb{R}$. By the proof of  Lemma
\ref{jordanproduct} it is immediate that whenever $(A_1, A_2,\ldots ,A_m)$ over $\mathbb{C}$ (respectively
$\mathbb{R}$) is hypercyclic one can completely describe the set of hypercyclic vectors as
\[
\left\{ \left(
\begin{array}{c}
z_1 \\
z_2
\end{array} \right) \in \mathbb{C}^2 : z_2\neq 0 \right\}\quad
\left(\textrm{respectively} \,\,\,
\left\{ \left(
\begin{array}{c}
x_1 \\
x_2
\end{array} \right) \in \mathbb{R}^2 : x_2\neq 0 \right\} . \right)
\]
\end{remark}

\subsection{The real case}
We shall need the following well known result, see for example \cite{Feldman}.
\begin{theorem} \label{Feldmandense}
If $a,b>1$ and $\frac{ \ln a}{\ln b}$ is irrational then the sequence $\{ \frac{a^n}{b^m}:n,m \in \mathbb{N} \}$
is dense in $\mathbb{R}^{+}$.
\end{theorem}

\begin{lemma}
Let $a,b \in\mathbb{R}$ such that $-1<a<0$, $b>1$ and $\frac{ \ln |a|}{\ln b}$ is irrational. Then the sequence
$\{ a^nb^m: n,m\in \mathbb{N} \}$ is dense in $\mathbb{R}$.
\end{lemma}
\begin{proof}
Since $\frac{ \ln |a|}{\ln b}$ is irrational it follows that $ \ln b /\ln \frac{1}{a^2}$ is irrational as well.
Applying Theorem \ref{Feldmandense} we conclude that the sequence $\{ a^{2n}b^m: n,m \in \mathbb{N} \}$ is dense
in $\mathbb{R}^{+}$. On the other hand the fact that $a$ is negative implies that the sequence $\{ a^{2n+1}b^m:
n,m \in \mathbb{N} \}$ is dense in $\mathbb{R}^{-}$. This completes the proof of the lemma.
\end{proof}
In what follows, for any $x\in{\mathbb R}$ we will be denoting by $[x]$ the `integral part of $x$', that is, the largest integer which does not exceed $x$, and
by $\{x\}$ the `fractional part of $x$', that is, $\{x\}=x-[x]$.

\begin{proposition}\label{R2dense}
There exist $a_1, a_2, a_3, a_4\in \mathbb{R} $ such that the sequence
\[ \left\{\left( \begin{array}{c}
\frac{k_1}{a_1}+\frac{k_2}{a_2}+\frac{k_3}{a_3}  +\frac{k_4}{a_4} \\
{a_1}^{k_1}{a_2}^{k_2}{a_3}^{k_3} {a_4}^{k_4}
\end{array} \right) : k_1, k_2, k_3 ,k_4 \in \mathbb{N}\right\} \]
is dense in $\mathbb{R}^2$.
\end{proposition}
\begin{proof}
By the lemma above fix $a,b\in\mathbb{R}$ such that $-1<a<0$, $a+\frac{1}{a}\in \mathbb{R}\setminus\mathbb{Q}$
and $\{ a^n b^m : n,m\in\mathbb{N}\}$ is dense in $\mathbb{R}$.  Let $x_1, x_2\in\mathbb{R}$ and $\epsilon>0$ be
given. Then there exist $n,m\in\mathbb{N}$ such that $|a^n b^m - x_2|<\epsilon$. Note that $a^n b^m= a^{n+k} b^m
\frac{1}{a^k} 1^s$ for every $k,s\in\mathbb{N}$. Note also that $a+\frac{1}{a}<0$. In the case $x_1\geq 0$
(resp. $x_1<0$) there exists $k\in\mathbb{N}$ such that $\frac{n}{a}+\frac{m}{b}+k\left (a+\frac{1}{a}\right )$
is less than $-1$ (resp. $x_1-1$) and
\[\left |\left\{\frac{n}{a}+\frac{m}{b}+k\left (a+\frac{1}{a}\right )\right\}-\{x_1\}\right |<\epsilon.\]
As there exists $s\in\mathbb{N}$ such that
\[\left|\frac{n}{a}+\frac{m}{b}+k\left (a+\frac{1}{a}\right )+s - x_1\right|<\epsilon \]
we are done.
Hence, setting $a_1=a, a_2=b, a_3=\frac{1}{a}, a_4=1$ the result is proved.
\end{proof}
\textit{Proof of Theorem \ref{R2hypjordan}}. This is an immediate consequence of Lemma \ref{jordanproduct},
Proposition \ref{R2dense} and Remark \ref{hypvectors}.

\begin{example}
One may construct many concrete examples of four $2\times 2$ matrices, in Jordan form over $\mathbb{R}$, being
hypercyclic. For example, fix $a,b \in\mathbb{R}$ such that $-1<a<0$, $b>1$ and both $a+\frac{1}{a}$, $\frac{
\ln |a|}{\ln b}$ are irrational. From the above we conclude that
\[ \left(
\left( \begin{array}{cc}
a & 1 \\
0 & a
\end{array} \right) ,
\left( \begin{array}{cc}
b & 1 \\
0 & b
\end{array} \right) ,
\left( \begin{array}{cc}
\frac{1}{a} & 1 \\
0 & \frac{1}{a}
\end{array} \right) ,
\left( \begin{array}{cc}
1 & 1 \\
0 & 1
\end{array} \right)
\right)
\]
is hypercyclic.
\end{example}

\begin{proposition} \label{2times2jordan}
(i) Every pair $(A_1, A_2)$ of $2\times 2$ matrices over $\mathbb{R}$ with $A_j$, $j=1,2$ being either diagonal
or in Jordan form is not hypercyclic.

(ii) There exist pairs $(A_1, A_2)$ of $2\times 2$ matrices over $\mathbb{R}$ such that $A_1$ is diagonal, $A_2$
is antisymmetric (rotation matrix) and $(A_1, A_2)$ is hypercyclic. In particular every non-zero vector in
$\mathbb{R}^2$ is hypercyclic for $(A_1, A_2)$, i.e.
\[ HC((A_1,A_2))=\mathbb{R}^2 \setminus \{ (0,0) \} .\]

(iii) There exist pairs $(A_1, A_2)$ of $2\times 2$ matrices over $\mathbb{R}$ such that both $A_1$ and $A_2$ are
antisymmetric and $(A_1, A_2)$ is hypercyclic. In particular every non-zero vector in $\mathbb{R}^2$ is
hypercyclic for $(A_1, A_2)$, i.e.
\[ HC((A_1,A_2))=\mathbb{R}^2 \setminus \{ (0,0) \} .\]

\end{proposition}
\begin{proof}
Let us prove assertion $(i)$. The case of $A_1, A_2$ both diagonal is covered by Feldman, see \cite{Feldman}. Assume that $A_1$ is diagonal and $A_2$ is in Jordan
form, i.e. ${A_1}=\left(
\begin{array}{cc}
a & 0\\
0 & a
\end{array} \right) ,\quad {A_2}=\left(
\begin{array}{cc}
b & 1\\
0 & b
\end{array} \right)$ for $a,b\in \mathbb{R}$.
Suppose that $(A_1,A_2)$ is hypercyclic and let $ \left( \begin{array}{c}
x_1 \\
x_2
\end{array} \right) \in \mathbb{R}^2$ be a hypercyclic vector for $(A_1,A_2)$. Then the sequence
\[\left\{ {A_1}^n{A_2}^m\left( \begin{array}{c}
x_1 \\
x_2
\end{array} \right) :n,m \in \mathbb{N} \right\},\]
i.e. the sequence $\left\{ \left( \begin{array}{c}
a^nb^mx_1+ma^nb^mx_2 \\
a^nb^mx_2
\end{array} \right) :n,m\in \mathbb{N} \right\}$
is dense in $\mathbb{R}^2$. Observe that $x_2$ cannot be
zero. Then there exist $y_1,y_2\in \mathbb{R}$ and sequences of positive integers $\{ n_k\}$, $\{ m_k\}$ such
that $m_k \to +\infty $ and \[ a^{n_k}b^{m_k}x_1+m_ka^{n_k}b^{m_k}x_2 \to y_1 ,\] \[ a^{n_k}b^{m_k}x_2 \to y_2\]
as $k\to +\infty$. It clearly follows that $|m_ka^{n_k}b^{m_k}x_2|\to +\infty $ which is a contradiction. Assume
now that both $A_1, A_2$ are in Jordan form, i.e.
\[
{A_1}=\left(
\begin{array}{cc}
a & 1\\
0 & a
\end{array} \right) ,\quad {A_2}=\left(
\begin{array}{cc}
b & 1\\
0 & b
\end{array} \right) ,\] for $a,b\in \mathbb{R}$ and $(A_1, A_2)$ is hypercyclic. Lemma \ref{jordanproduct} implies
that the sequence \[ \left\{\left( \begin{array}{c}
\frac{n}{a}+\frac{m}{b} \\
a^nb^m
\end{array} \right) : n,m \in \mathbb{N}\right\} \]
is dense in $\mathbb{R}^2$. Observe that both $|a|,|b|$ are not equal to $1$. By taking absolute value in the
second coordinate and then applying the logarithmic function, it follows that the sequence
\[ \left\{\left( \begin{array}{c}
\frac{n}{a}+\frac{m}{b} \\
n\ln |a|+m\ln |b|
\end{array} \right) : n,m \in \mathbb{N}\right\} \]
is dense in $\mathbb{R}^2$. Hence the sequence
\[ \left\{\left( \begin{array}{c}
n\frac{\ln |a|}{a}+m\frac{\ln |a|}{b} \\
n\frac{\ln |a|}{a}+m\frac{\ln |b|}{a}
\end{array} \right) : n,m \in \mathbb{N}\right\} \]
is dense in $\mathbb{R}^2$. Subtracting the second coordinate from the first one we conclude that the sequence
\[ \left\{ m\left( \frac{\ln |a|}{b}-\frac{\ln |b|}{a}\right) : m \in \mathbb{N} \right\} \] is dense
in $\mathbb{R}$ which is absurd. We proceed with the proof of assertion $(ii)$. There exist $a\in
\mathbb{R}\setminus \mathbb{Q}$ and $b\in \mathbb{C}$ such that the sequence $\{ a^nb^m: n,m \in \mathbb{N} \}$
is dense in $\mathbb{C}$, see \cite{Feldman}. Write $b=|b|e^{i\theta }$ and set
\[
{A_1}=\left(
\begin{array}{cc}
a & 0\\
0 & a
\end{array} \right) ,\quad {A_2}=\left(
\begin{array}{cc}
|b|\cos \theta & -|b|\sin \theta\\
 |b|\sin \theta & |b|\cos \theta
\end{array} \right) .\]
Then we have
\[
{A_1}^n{A_2}^m=\left(
\begin{array}{cc}
a^n|b|^m \cos (m\theta )& -a^n|b|^m \sin (m\theta )\\
a^n|b|^m \sin (m\theta ) & a^n|b|^m \cos (m\theta )
\end{array} \right) .\]
Applying in the above relation the vector $ \left( \begin{array}{c}
1 \\
0
\end{array} \right) $ and taking into account that the sequence $\{ a^nb^m: n,m \in \mathbb{N} \}$
is dense in $\mathbb{C}$ we conclude  that the sequence
\[ \left\{
{A_1}^n{A_2}^m \left( \begin{array}{c}
1 \\
0
\end{array} \right) :n,m\in \mathbb{N}  \right\}= \left\{   \left(
\begin{array}{c}
a^n|b|^m \cos (m\theta ) \\
a^n|b|^m \sin (m\theta )
\end{array} \right) :n,m\in \mathbb{N} \right\} \]
is dense in $\mathbb{R}^2$. Hence $(A_1, A_2)$ is hypercyclic. It is now easy to show that every non-zero vector
in $\mathbb{R}^2$ is hypercyclic for $(A_1,A_2)$.

In order to prove the last assertion we follow a similar line of reasoning as above. That is, by Corollary 3.2
in \cite{Feldman} there exist $a,b\in \mathbb{C} \setminus \mathbb{R}$ such that the sequence $\{ a^nb^m: n,m
\in \mathbb{N} \}$ is dense in $\mathbb{C}$. Write $a=|a| e^{i\phi} $, $b=|b|e^{i\theta }$ and set
\[
{A_1}=\left(
\begin{array}{cc}
|a|\cos \phi & -|a|\sin \phi\\
|a|\sin \phi & |a|\cos \phi
\end{array} \right) ,\quad {A_2}=\left(
\begin{array}{cc}
|b|\cos \theta & -|b|\sin \theta\\
 |b|\sin \theta & |b|\cos \theta
\end{array} \right) .\] A direct computation gives that
$\left\{ {A_1}^n{A_2}^m \left( \begin{array}{c}
1 \\
0
\end{array} \right) :n,m\in \mathbb{N}  \right\}$ equals to $\left\{   \left(
\begin{array}{c}
|a|^n|b|^m \cos (n\phi -m\theta ) \\
|a|^n|b|^m \sin (n\phi +m\theta )
\end{array} \right) :n,m\in \mathbb{N} \right\}$
and by the choice of $a,b$ we conclude that the vector $ \left( \begin{array}{c}
1 \\
0
\end{array} \right) $ is hypercyclic for $(A_1,A_2)$. This completes the proof of the proposition.
\end{proof}

\begin{question}
What is the minimum number of $2\times 2$ matrices over $\mathbb{R}$ in Jordan form so that their tuple form a
hypercyclic operator?
\end{question}

\subsection{The complex case}
In what follows we will be writing $\Re(z)$
and $\Im(z)$ for the real and imaginary parts of a complex number $z$ respectively.
\begin{proposition} \label{complexJordan}
There exist $a_j\in \mathbb{C} $, $j=1,2,\ldots ,8$  such that the sequence
\[ \left\{\left( \begin{array}{c}
\frac{k_1}{a_1}+\frac{k_2}{a_2}+\ldots  +\frac{k_8}{a_8} \\
{a_1}^{k_1}{a_2}^{k_2}\ldots {a_8}^{k_8}
\end{array} \right) : k_1, k_2, \ldots ,k_8 \in \mathbb{N}\right\} \]
is dense in $\mathbb{C}^2$.
\end{proposition}
\begin{proof}
The proof is in the same spirit as the proof of Proposition \ref{R2dense}. Fix $a,b\in\mathbb{C}$ such that
$-1<a<0$, $a+\frac{1}{a},a-\frac{1}{a}\in \mathbb{R}\setminus\mathbb{Q}$ and $\{ a^n b^m : n,m\in\mathbb{N}\}$
is dense in $\mathbb{C}$ (see Corollary \ref{denseinCR}). Let $z_1, z_2\in\mathbb{C}$ and $\epsilon>0$ be given.
Then there exist $n,m\in\mathbb{N}$ such that $|a^n b^m - z_2|<\epsilon$. Note that
\[a^n b^m= a^{n+k} b^m \frac{1}{a^k} 1^s (ia)^{\xi}\left(\frac{1}{ia}\right)^{\xi}
(4i)^{\rho}\left(-\frac{1}{4}\right)^{\rho}\] for every $k,s,\xi\in\mathbb{N}$ and $\rho\in 4{\mathbb N}$.
Note that $a+\frac{1}{a}<0$ and $a-\frac{1}{a}>0$. In the case $\Re(z_1)\geq 0$ (resp. $\Re(z_1)<0$)
there exists $k\in\mathbb{N}$ such that
$\Re\left(\frac{n}{a}+\frac{m}{b}+k\left(a+\frac{1}{a}\right )\right)$ is less than $-1$ (resp. $\Re(z_1)-1$) and
\[\left |\left\{\Re\left(\frac{n}{a}+\frac{m}{b}+k\left (a+\frac{1}{a}\right )\right)\right\}-\{\Re(z_1)\}\right |<\epsilon.\]
Also, in the case $\Im(z_1)\geq 0$ (resp. $\Im(z_1)<0$) there exists $\xi\in{\mathbb N}$ such that
$\Im\left(\frac{n}{a}+\frac{m}{b}+i\xi\left(a-\frac{1}{a}\right )\right)$ is greater than $1$ (resp. $\Im(z_1)+1$) and
\[\left |\left\{\Im\left(\frac{n}{a}+\frac{m}{b}+i\xi\left (a-\frac{1}{a}\right )\right)\right\}-\{\Im(z_1)\}\right |<\epsilon.\]
As there exists $\rho\in 4\mathbb{N}$ such that
\[\left|\Im\left(\frac{n}{a}+\frac{m}{b}+i\xi\left (a-\frac{1}{a}\right )-\rho\left(\frac{i}{4}\right)\right) - \Im(z_1)\right|<\epsilon \]
we are done with the imaginary part. Also, there exists $s\in{\mathbb N}$ such that
\[\left|\Re\left(\frac{n}{a}+\frac{m}{b}+k\left(a+\frac{1}{a}\right)-4\rho+s\right) - \Re(z_1)\right|<\epsilon. \]
But this means that the real and imaginary parts of the complex number
\[\frac{n}{a}+\frac{m}{b}+k\left(a+\frac{1}{a}\right)+s+i\xi\left(a-\frac{1}{a}\right)-\rho\frac{i}{4}-4\rho\]
are within $\epsilon$ of the real and imaginary parts of $z_1$.
Hence, setting $a_1=a, a_2=b, a_3=\frac{1}{a}, a_4=1, a_5=ia, a_6=\frac{1}{ia}, a_7=4i, a_8=-\frac{1}{4}$
the result is proved.
\end{proof}
\textit{Proof of Theorem \ref{C2hypjordan}}. By Proposition \ref{complexJordan}, Lemma \ref{jordanproduct} and
Remark \ref{hypvectors} the assertion follows.

\begin{example}
Fix $a,b\in\mathbb{C}$ such that $-1<a<0$, $a+\frac{1}{a},a-\frac{1}{a}\in \mathbb{R}\setminus\mathbb{Q}$ and
$\{ a^n b^m : n,m\in\mathbb{N}\}$ is dense in $\mathbb{C}$. From the above it is evident that the following
$8$-tuple of $2\times 2$ matrices in Jordan form over $\mathbb{C}$
\[
\left( \begin{array}{cc}
a & 1 \\
0 & a
\end{array} \right) ,
\left( \begin{array}{cc}
b & 1 \\
0 & b
\end{array} \right) ,
\left( \begin{array}{cc}
\frac{1}{a} & 1 \\
0 & \frac{1}{a}
\end{array} \right) ,
\left( \begin{array}{cc}
1 & 1 \\
0 & 1
\end{array} \right) ,
\]
\[
\left( \begin{array}{cc}
ia & 1 \\
0 & ia
\end{array} \right) ,
\left( \begin{array}{cc}
\frac{1}{ia} & 1 \\
0 & \frac{1}{ia}
\end{array} \right) ,
\left( \begin{array}{cc}
4i & 1 \\
0 & 4i
\end{array} \right) ,
\left( \begin{array}{cc}
-\frac{1}{4} & 1 \\
0 & -\frac{1}{4}
\end{array} \right)
\]
is hypercyclic.

\end{example}

\begin{question}
What is the minimum number of $2\times 2$ matrices over $\mathbb{C}$ in Jordan form so that their tuple form a
hypercyclic operator?
\end{question}

\section{Products of $3\times 3$ matrices}
\begin{proposition} {\bf (Feldman)} \label{Feldman}
If $b_1,b_2 \in \mathbb{D} \setminus \{ 0\}$ then there exists a dense set $\Delta \subset \mathbb{C} \setminus
\mathbb{D} $ such that for every $a_1,a_2 \in \Delta $ the sequence
\[ \left\{\left( \begin{array}{c}
{a_1}^n{b_1}^m \\
{a_2}^n{b_2}^l
\end{array} \right) : n,m,l \in \mathbb{N}\right\} \]
is dense in $\mathbb{C}^2$.
\end{proposition}

\begin{corollary}\label{denseinCR}
There exist $a\in \mathbb{C}$ and $b,c,d \in \mathbb{R}$ such that the sequence
\[ \left\{\left( \begin{array}{c}
a^nb^m \\
c^nd^l
\end{array} \right) : n,m,l \in \mathbb{N}\right\} \]
is dense in $\mathbb{C} \times \mathbb{R}$.
\end{corollary}
\begin{proof}
Fix two real numbers $b_1, b_2$ with $b_1,b_2 \in (0,1)$. By Proposition \ref{Feldman} there exist $a_1, a_2\in
\mathbb{C} \setminus \mathbb{D}$ such that the sequence
\[ \left\{\left( \begin{array}{c}
{a_1}^n{b_1}^m \\
{a_2}^n{b_2}^l
\end{array} \right) : n,m,l \in \mathbb{N}\right\} \]
is dense in $\mathbb{C}^2$. Define $a=a_1$, $b=b_1$, $c=|a_2|$ and $d=-\sqrt{b_2}$. Observe that the sequence
\[ \left\{\left( \begin{array}{c}
a^n b^m \\
c^n {b_2}^l
\end{array} \right) : n,m,l \in \mathbb{N}\right\} \]
is dense in $\mathbb{C} \times [0, +\infty )$. Take $z\in \mathbb{C}$ and $x\in \mathbb{R}$.\\
{\bf Case I}: $x \geq 0$.\\
Then there exist sequences of positive integers $\{ n_k\} ,\{ m_k \}$, $\{ l_k\}$ such that \[ a^{n_k}b^{m_k}\to
z \,\,\, \textrm{and} \,\,\, c^{n_k}{b_2}^{l_k}\to x.\] Since ${b_2}^{l_k}=d^{2l_k}$ we get $c^{n_k}d^{2l_k}\to
x$.\\
{\bf Case II}: $x <0$.\\
Then there exist sequences of positive integers $\{ n_k\} ,\{ m_k \}$, $\{ l_k\}$ such that \[ a^{n_k}b^{m_k}\to
z \,\,\, \textrm{and} \,\,\, c^{n_k}{b_2}^{l_k}\to \frac{x}{d}.\] The last implies that $c^{n_k}d^{2l_k+1}\to
x$. This completes the proof of the corollary.
\end{proof}

\begin{proposition}\label{3times3real}
There exist $3$ tuples $(A_1, A_2, A_3)$ of $3\times 3$ matrices over $\mathbb{R}$ such that $(A_1, A_2, A_3)$
is hypercyclic.
\end{proposition}
\begin{proof}
By Corollary \ref{denseinCR} there exist $a\in \mathbb{C}$ and $b,c,d \in \mathbb{R}$ such that the sequence
\[ \left\{\left( \begin{array}{c}
a^nb^m \\
c^nd^l
\end{array} \right) : n,m,l \in \mathbb{N}\right\} \]
is dense in $\mathbb{C} \times \mathbb{R}$. Write $a=|a|e^{i\theta }$ and set
\[
{A_1}=\left(
\begin{array}{ccc}
|a|\cos \theta & -|a|\sin \theta & 0\\
|a|\sin \theta & |a|\cos \theta & 0 \\
0 & 0& c
\end{array} \right) ,\,\, {A_2}=\left(
\begin{array}{ccc}
b & 0 & 0\\
0 & b & 0\\
0 & 0 & 1
\end{array} \right) \,\,\mbox{and}\]
${A_3}=\left(
\begin{array}{ccc}
1 & 0 & 0\\
0 & 1 & 0\\
0 & 0 & d
\end{array} \right).$
Then we have
\[ {A_1}^n{A_2}^m{A_3}^l= \left( \begin{array}{ccc}
|a|^nb^m\cos (n\theta) & -|a|^nb^m\sin (n\theta) & 0\\
|a|^nb^m\sin (n\theta) & |a|^nb^m\cos (n\theta) & 0 \\
0 & 0& c^nd^l
\end{array} \right) ,\]
which in turn gives
\[ {A_1}^n{A_2}^m{A_3}^l \left( \begin{array}{c}
1 \\
0 \\
1
\end{array} \right) = \left( \begin{array}{c}
|a|^nb^m\cos (n\theta) \\
|a|^nb^m\sin (n\theta) \\
c^nd^l \end{array} \right) .\] The last and the choice of $a,b,c,d$ imply that $(A_1,A_2,A_3)$ is hypercyclic
with $\left( \begin{array}{c}
1 \\
0 \\
1
\end{array} \right)$ being a hypercyclic vector for $(A_1,A_2,A_3)$.
\end{proof}

\section{Proof of Theorem \ref{maintheorem}}
By Proposition \ref{2times2jordan}, there exist $2\times 2$ matrices $B_1$ and $B_2$ such that $(B_1, B_2)$ is
hypercyclic.\\
{\bf Case I}: $n=2k$ for some positive integer $k$. For $k=1$ the result follows by Proposition
\ref{2times2jordan}. Assume that $k>1$. Each $A_j$ will be constructed by blocks of $2\times 2$ matrices. Let
$I_2$ be the $2\times 2$ identity matrix. We will be using the notation $diag(D_1, D_2,\cdots,D_n)$
to denote the diagonal matrix with diagonal entries the block matrices $D_1, D_2,\cdots, D_n$.
Define $A_1=diag(B_1, I_2,\cdots, I_2)$, $A_2=diag(B_2, I_2,\cdots, I_2)$, $A_3=diag(I_2, B_1, I_2,\cdots, I_2)$,
$A_4=diag(I_2, B_2, I_2,\cdots,I_2)$ and so on up to $A_{n-1}=diag(I_2,\cdots, I_2, B_1)$,
$A_n=diag(I_2,\cdots, I_2, B_2)$.
It is now easy to check that $(A_1, A_2, \ldots A_n)$ is hypercyclic and furthermore that the set
$HC((A_1,A_2, \ldots , A_n))$ is \[ \{ (x_1, x_2, \ldots ,x_n) \in \mathbb{R}^n: x_{2j-1}^2+x_{2j}^2\neq 0 ,\forall
j=1,2,\ldots ,k \} .\]
{\bf Case II}: $n=2k+1$ for some positive integer $k$. If $k=1$ the result follows by Proposition
\ref{3times3real}. Suppose $k>1$. For simplicity we treat the case $k=2$, since the general case follows by
similar arguments. By Proposition \ref{3times3real} there exist $C_1, C_2, C_3$, $3\times 3$ matrices such that
$(C_1, C_2, C_3)$ is hypercyclic. Let $I_3$ be the $3\times 3$ identity matrix. Define
$A_1=diag(B_1, I_3)$, $A_2=diag(B_2, I_3)$, $A_3=diag(I_2, C_1)$, $A_4=diag(I_2, C_2)$ and $A_5=diag(I_2, C_3)$.
It can easily be shown that $( A_1, A_2, \ldots ,A_5 )$ is hypercyclic. The details are left to the reader.

\end{document}